\documentstyle[12pt,amssymb,pb-diagram,CD,papersaver6]{article}%12pt,amscd

\def\afooter#1#2{{\def\thefootnote{\mbox{${}^{#1}$}}\mbox{\footnotemark[0]}
        \footnotetext{#2} }}
\newtheorem{Theorem}{Theorem}[section]
\newtheorem{Definition}{Definition}[section]

\newtheorem{Lemma}[Theorem]{Lemma}
\newtheorem{Proposition}[Theorem]{Proposition}

\newtheorem{Corollary}[Theorem]{Corollary}

\newtheorem{Remark}{Remark}[section]
\newtheorem{fact}[Theorem]{Fact}

\let\qed=\Box
\def\){ \right) }
\def\({ \left( }
\def\[{ \left[ }
\def\]{ \right] }
\def\<{ \langle }
\def\>{ \rangle }
\let\ljunk=\{
\let\rjunk=\}
\def\{{\left\ljunk}
\def\}{\right\rjunk}
\def\p{\partial}

\def\Diff{\mbox{Diff}}

\def\mod{\mbox{mod }}

\newcommand{\Z}{{\mathbf Z}}

\newcommand{\Q}{{\mathbf Q}}
\newcommand{\Cat}{{\mathrm C}{\mathrm a}{\mathrm t}}
\newcommand{\PL}{{\mathrm P}{\mathrm L}}
\newcommand{\Top}{{\mathrm T}{\mathrm o}{\mathrm p}}
\newcommand\halfp{\frac{p_1}{2}}

\newcommand\NO{{\mathcal{N}}^{O}}
\newcommand\SO{{\mathcal{S}}^{O}}

\newcommand{\h}{{\mbox{Hom}}}
\begin{document} 

\title{A classification of $S^3$-bundles over
$S^4$\afooter{\mbox{ \ }}{Keywords and phrases. sphere bundles over
spheres,
 $\mu$-invariant, homotopy, homeomorphism,
       diffeomorphism classification\\ \hspace*{5mm} 2000 Mathematics
Subject 
Classification.  55R15,55R40,57T35}} \author{Diarmuid Crowley, 
Christine M. Escher}
\date{April 22, 2000}

\maketitle

\begin{abstract}
We classify the total spaces of bundles over the four 
sphere with fiber a three sphere up to orientation preserving and
reversing
homotopy equivalence, 
	homeomorphism and diffeomorphism.  These total spaces have
	been of interest to both topologists and geometers.  It has
recently 
	been shown by Grove and Ziller \cite{GZ} that each of these total
spaces
	admits metrics with nonnegative sectional curvature.
\end{abstract}

\section{Introduction}
For almost fifty years, $3$-sphere bundles over the $4$-sphere have been
of 
interest to both topologists and geometers.  In 1956, Milnor \cite{M} 
proved that all  $S^3$-bundles over $S^4$ with Euler class $e = \pm 
1$ are homeomorphic to $S^{7}$.  He also showed that some of these 
bundles are not diffeomorphic to $S^{7}$ and thereby exhibited the first 
examples of exotic spheres.  Shortly thereafter, in 1962, Eells and Kuiper 
\cite{EK} classified all such bundles with Euler class $e = \pm 
1$ up to diffeomorphism and showed that $15$ of the $27$ seven 
dimensional exotic spheres can be described as $S^3$-bundles over $S^4$.

In 1974 Gromoll and Meyer \cite{GM} constructed a metric with 
non-negative sectional curvature on one of these sphere bundles, 
exhibiting the first example of an exotic sphere with non-negative 
curvature.  Very recently Grove and Ziller \cite{GZ} showed that the 
total space of every $S^3$-bundle over $S^4$ admits a metric with 
non-negative sectional curvature.  Motivated by these examples they 
asked for a classification of these manifolds up to homotopy 
equivalence, homeomorphism and diffeomorphism (see problem 5.2 in 
\cite{GZ}).  Although partial classifications in various
categories appear in \cite{JW2}, \cite{S}, \cite{T} and
\cite{Wi}, it seems that the complete classification has not been carried
out.  
It is the  purpose of this article to do this.

We consider fiber 
bundles over the $4$-sphere $S^{4}$ with total space $M$, fiber the 
$3$-sphere $S^{3}$ and structure group $SO(4)$.  Equivalence 
classes of such bundles are in one-to-one correspondence with 
$\pi_{3}(SO(4)) \cong {\mathbf Z} \oplus {\mathbf Z}$.  Following 
James and Whitehead \cite{JW2}, we choose generators $\rho$ and $\sigma$ 
of $\pi_{3}(SO(4))$ such that
$$
\rho(u) \ v = u \ v \ u^{-1} \ , \quad \sigma(u) \ v = u \ v \ ;
$$
where $u$ and $v$ denote the quaternions with norm 1, i.e. we have
have identified $S^{3}$ with the unit quaternions.  With this choice
of generators a pair of integers $(m,n)$ gives the element 
$m \ \rho + n \ \sigma
\in \pi_{3}(SO(4))$  and thus determines both a vector bundle 
$\xi_{m,n}$ and the corresponding sphere bundle 
$\pi_{m,n}: M_{m,n}:=S(\xi_{m,n}) \longrightarrow S^{4}$.

\begin{Remark} The definition we use is different from the one
given by Milnor in {\rm \cite{M}}. He uses two integers $(k,l)$
corresponding to a different choice of generators of $ \pi_{3}(SO(4))
\cong \Z\oplus \Z$. The two pairs $(m,n)$ and $(k,l)$ are related by
$k+l=n$, $l=-m$.  

Note that a change of orientation in the fiber 
leads to $M_{m,n} \cong M_{m+n,-n}$ whereas a change of orientation 
in the base leads to $M_{m,n} \cong M_{-m,-n}$.  Hence we now always 
assume that $n \ge 0$. 
We can also exclude the case of $n=0$ as $M_{m,0}$ is diffeomorphic or 
homeomorphic to $M_{m',0}$ only when $m'=\pm m$.   This follows from
the 
topological invariance of the rational Pontrjagin
classes and the fact that $M_{m,n}$ is diffeomorphic to $M_{-m-n,n}$ 
for any $m$ and $n$.  James and
Whitehead \cite{JW1} proved that $M_{m,0}$ is homotopic to $M_{m',0}$ if
and only if
$m' \equiv \pm m ~\mod 12$.  Lastly, each $M_{m,0}$ has an orientation
reversing 
self-diffeomorphism. Hence we now always assume that $n > 0$.
\end{Remark}

The first level of classification of the manifolds $M_{m,n}$ is up
to homotopy equivalence. In a remarkable paper \cite{JW2} James and
Whitehead succeeded in  classifying the pairs of
manifolds $(M_{m,n},S^3)$ up to homotopy equivalence except for the case
$n=2$.  A few years 
later Sasao \cite{S} undertook the homotopy classification of the total
spaces $M_{m,n}$.  Let $M(Z_n,3)$ denote the Moore space formed by
attaching a $4$-disc to a
$3$-sphere with a degree $n$ map of the $3$-sphere.  Each $M_{m,n}$ has
the
homotopy type of $M(Z_n,3)$ with a $7$-cell attached. Sasao
computed both $\pi_6(M(Z_n,3))$ and the action of $[M(\Z_n,3),M(\Z_n,3)]$
on
$\pi_6(M(\Z_n,3))$.  In principle, this should allow for a complete
homotopy
classification though Sasao did not do the computation.  We deduce the
homotopy classification as a consequence of the homeomorphism
classification
and thus avoid the homotopy theory completely.

\begin{Theorem}\label{hom}

Let $n , n'> 0$.

\begin{enumerate}
\item  
The manifolds  $M_{m',n'}$ and $M_{m,n}$
are  orientation preserving  homotopy equivalent if and only if $n = n'$
and 
$m' \equiv \alpha \,m ~ \mod (n,12)$ where $\alpha^2 \equiv 1 ~ \mod 
(n, 12)$.

\item Orientation reversing homotopy equivalences 
between any $M_{m',n }$ and $M_{m,n}$ can only exist when
$n=2^{\epsilon}p_1^{i_1}
\dots p_k^{i_k},~ p_{i} ~ \mbox{prime}, ~p_i \equiv 1 ~ \mod 4$ 
and $\epsilon = 0,1.$  
Furthermore if $n$ is
of this form with 
$\epsilon = 0$, then the single oriented homotopy type admits an
orientation
reversing self homotopy equivalence; if $\epsilon = 1,$ $M_{m',n }$ 
is orientation reversing homotopy equivalent to $M_{m,n}$ if and only if
$m'+m \not\equiv 0~\mod 2.$  

\end{enumerate}
\end{Theorem}

Here $(n,12)$ denotes the greatest common divisor of $n$ and $12$. 
Notice that since $H^4(M_{m,n};\Z)\cong \Z_{n}$ and we assumed $n$ to be 
nonnegative, necessarily $n=n^{\prime}$ if $M_{m,n}$ and
$M_{m^{\prime},n^{\prime}}$ are homotopy equivalent.

\begin{Remark}  Comparing Theorem \ref{hom} with \cite{JW2} (for $n > 
2$, $(M_{m,n},S^3)$ is homotopy equivalent to $(M_{m',n},S^3)$ if and 
only if $m' \equiv \,\pm \,m ~ \mod (n,12)$) one 
concludes that the homotopy classification of the manifold pairs
$(M_{m,n},S^3)$  differs from the homotopy classification of the total
spaces
$M_{m,n}$ if and only if $n$ is divisible by $12$.  In this case the 
congruence $\alpha^2 \equiv 1 ~ \mod (n, 12)$ has two nontrivial 
solutions, namely $\alpha \equiv \pm 5 ~ \mod 12\,.$
\end{Remark}

The next level, the homeomorphism classification, was first studied by
Tamura \cite{T}, who constructed explicit homeomorphisms between $M_{m,n}$
and
$M_{m',n}$ if $m \equiv \pm m'~\mod n$.  These homeomorphisms were
constructed with the aid of specific foliations of the total spaces.
However, as handlebody theory developed, general techniques for
classifying 
highly connected manifolds became available and Wall (\cite{Wa1}) was able
to classify $(s-1)$ connected $(2s+1)$ manifolds except when $s=3,7$.
Wilkens \cite{Wi} then extended the techniques of Wall to complete the
cases
$s=3,7$ except for an occasional $\Z_{2}$-ambiguity.  We complete Wilkens'
classification in the case of $S^3$-bundles over $S^4$ using two
topological
invariants.  The first invariant was used by Wilkens and is the
characteristic
class of Spin manifolds, $\halfp$, which is defined as the
generator
of $H^4(Bspin;\Z)$ that has the same sign as the first Pontrjagin 
class $p_1$.  We stress that the
division
by two occurs universally in $H^4(Bspin;\Z)$ so that $\halfp(M_{m,n})$
should not
be thought of as $p_1(M_{m,n})/2$ which does not make sense in the
presence
of $2$-torsion.  Pulling back a preferred generator of
$H^4(S^4;\Z)$ to $H^4(M_{m,n};\Z)$ yields a canonical identification of
$H^4(M_{m,n};\Z)$ with $\Z_n$ and hence we identify $\halfp(M_{m,n})$ with
an
element of $\Z_n$.  The second invariant is the topological Eells-Kuiper
invariant $s_{1}(M_{m,m}) := 28 \mu(M_{m,n}) \in \Q/\Z$ which was used
in \cite{KS} where it is shown (Theorem 2.5) to be an invariant of
topological spin manifolds.  We now recall the definition of the $\mu$ 
invariant from \cite{EK}.

Assume that $M = M^{7}$ is a closed, smooth spin manifold such that
$M$ is a spin boundary, i.e.  $M =\p W^8$, where $W^{8}$ is a 
closed, smooth spin manifold. In addition, we require $(W,M)$ to
satisfy the {\it $\mu$-condition} which is automatic when $M$ is a
rational homology sphere, see \cite{EK}. The
$\mu$-condition allows the pullback of the Pontrjagin classes of $W$ to
$H^{*}(W,M)$. Then the $\mu$-invariant is defined as
\begin{equation}\label{Inv}
\mu(M^{7}) \equiv \frac{1}{2^{7} \cdot 7} \ \{p_{1}^{2}(W) - 4 \ 
\sigma[W]\} ~ \mod 1  
\end{equation}
Here $\sigma[W]$ stands for the signature of $W$. The most important feature
of
the $\mu$-invariant is that it is an invariant of the diffeomorphism
type of $M$ whereas $s_1:=28\mu$ is an invariant of the homeomorphism type
of $M$, see \cite{KS}.  Calculating these invariants for the total spaces
$M_{m,n}$ 
we obtain
$$\halfp(M_{m,n}) \equiv 2\,m ~ \mod n ; ~ s_{1}(M_{m,n}) \equiv
\frac{1}{8 n} (4\,m\,(m + n) + n \,(n-1))  ~ \mod 1 \,.$$

\begin{Theorem}\label{top} 
The following are equivalent.
\begin{enumerate}
\item  $M_{m',n}$ is orientation preserving [reversing] 
       $PL$-homeomorphic to $M_{m,n}.$
\item  $M_{m',n}$ is orientation preserving [reversing] homeomorphic to 
       $M_{m,n}.$
\item  $\halfp(M_{m',n}) = \alpha \, \halfp(M_{m,n})$ and $s_1(M_{m',n}) =
s_1(M_{m,n})$ where $\alpha^2 \equiv \,1 ~ \mod n$.

 [$\halfp (M_{m',n}) = \alpha ~ \halfp (M_{m,n})$ and $s_1(M_{m',n})
= - s_1(M_{m,n})$ where $\alpha^2 \equiv \, - 1 ~ \mod n$.]
\end{enumerate}
\end{Theorem}

Combining the congruences we obtain the following version of Theorem 
\ref{top}.

\begin{Corollary}
\label{topcor}
Let $n=2^a \, q,~q~$ odd.  
\begin{enumerate}
\item  $M_{m',n}$ is
orientation preserving homeomorphic to $M_{m,n}$ if and only if the
following
condition holds.
\begin{itemize}
\item $a = 0$: $m' \equiv \alpha \, m ~ \mod n $ where  $\alpha^2 \equiv
\,1 ~ \mod n$. 
\item $a=2$ and  $m$ odd, or $a > 2$ and $m$ even: 
$m' \equiv \alpha \, m ~ \mod \frac{n}{2}$ where  $\alpha^2 \equiv
\,1 ~ \mod n$. 
\item $a = 1$, or $a=2$ and $m$ even, or $a > 2$ and $m$ odd: 
$m' \equiv \alpha \, m ~\mod n$ where $\alpha \equiv \pm \,1 ~\mod
2^a$ 
and  $\alpha^2 \equiv \,1\, ~ \mod n$. 
\end{itemize}
\item  $M_{m',n}$ is
orientation reversing homeomorphic to $M_{m,n}$ if and only if the
following
condition holds.
\begin{itemize}
\item $n=p_1^{i_1} \dots p_k^{i_k},~p_i \equiv 1 ~ \mod 4$ and
  $m' \equiv \alpha \, m ~ \mod n $ where  $\alpha^2 \equiv
\,- 1 ~ \mod n ~$ or  
\item $n=2p_1^{i_1} \dots p_k^{i_k},~p_i \equiv 1 ~ \mod 4$ and 
 $m' \equiv \alpha \, (m + \frac{n}{2}) ~ \mod n $ where  $\alpha^2
\equiv
\,- 1 ~ \mod n$. 
\end{itemize}
\end{enumerate}
\end{Corollary}

The following corollary of the proof of Theorem \ref{top} generalizes 
a result of Kitchloo and Shankar \cite{KiSh}.

\begin{Corollary}
\label{topcor2}
Any $2$-connected $7$-manifold $M$ with $H^4(M;\Z) \cong \Z_n$ and
linking form isomorphic to $\,l: \Z_{n} \times \Z_{n} \longrightarrow 
\Q/\Z, (r,s) \mapsto \frac{rs}{n}\,$ is homeomorphic to the total space
of a $3$-sphere bundle over the $4$-sphere.
\end{Corollary}

The diffeomorphism classification of the manifolds $M_{m,n}$ 
follows the pattern that Eells and Kuiper discovered in the case $n=1$ 
for Milnor's exotic spheres.  That is, $M_{m',n}$ and $M_{m,n}$ are
diffeomorphic if and only if they are homeomorphic and their $\mu$
invariants coincide.

\begin{Theorem}\label{diff}
	The total space $M_{m',n}$ is orientation preserving
[reversing] diffeomorphic to 
 $M_{m,n}$ if and only if the following conditions hold.
 
 \medskip
 
 \noindent
 $\halfp(M_{m',n}) = \alpha \, \halfp(M_{m,n})$ and
$\mu(M_{m',n}) =
\mu(M_{m,n})$ where $\alpha^2 \equiv \,1 ~ \mod n$.
	   
\medskip 

\noindent
[$\halfp (M_{m',n}) = \alpha ~ \halfp (M_{m,n})$ and
$\mu(M_{m',n}) = -\mu(M_{m,n})$ where $\alpha^2 \equiv \, -
1 ~ \mod n$.]
\end{Theorem}

\medskip 

Using the definition of $\halfp$ and $\mu$ we describe Theorem 
\ref{diff} in terms of congruences as follows.

\begin{Corollary}
The total space $M_{m',n}$ is orientation preserving
[reversing] diffeomorphic to 
 $M_{m,n}$ if and only if the following conditions hold.

 \medskip
 
 \noindent
 $m'\,(n+m') \equiv m\,(n+m)  \, \mod 56\, n\,$ and 
 $2\,m' \equiv 2\,\alpha \,m ~ \mod n$ where 
       $\alpha^2 \equiv \,1 \, \mod n\,.$
	   
\medskip 

\noindent
 [$m'\,(n+m') \equiv - m\,(n+m) \, \mod 56 \,n\,$ and 
 $2\,m' \equiv 2\,\alpha \,m ~ \mod n$ where 
$\alpha^2 \equiv \,- 1 \, \mod n\,.$]
\end{Corollary}

Smooth surgery theory (see \cite{MM}) implies that there are exactly
28 different smooth manifolds homeomorphic to $M_{m,n}$. A
natural question to ask is which manifolds homeomorphic to
$M_{m,n}$ may be represented by manifolds $M_{m^{\prime},n}$.  In section 4
we answer this question for special cases of $n$.  We see that in 
some cases all $28$ differentiable structures are realized by 
$S^{3}$-bundles over $S^{4}$.  We also 
discuss the case $n=10$ which includes the seven dimensional Berger 
manifolds as shown in \cite{KiSh}.

\medskip

The first author wishes to thank Jim Davis for his patient guidance and
Fred Wilhelm for insightful conversations.  The second author would like
to thank Wolfgang Ziller for many helpful discussions.  Both authors are
indebted to Boris Botvinnik for his strong support and generosity.

\section{Outline of proof}

\subsection{A note on orientation}

Whereas classifications arising from the use of surgery theory
are naturally given up to 
orientation preserving maps, geometers often require 
classifications irrespective of orientation.
For our sphere bundles there are diffeomorphisms that reverse the 
orientation: of the fiber; of the base; or of both the fiber and the 
base.  These diffeomorphisms give rise to orientation reversing 
diffeomorphisms from $M_{m,n}$ to $M_{m+n,-n}$ and from $M_{m,n}$ to 
$M_{m,-n}$ respectively, as well as 
orientation preserving diffeomorphisms between $M_{m,n}$ and 
$M_{-m-n,n}$.  Hence it is sufficient to work with $n\geq0$; this is 
the convention established by James and Whitehead, \cite{JW2}.  
We now fix a generator 
$\iota_4$ of $H^4(S^4;\Z)$ and use the same notation for the images 
of $\iota_4$ under the isomorphisms 
$H^4(S^4;\Z) \cong H_4(S^4;\Z) \cong \pi_4(S^4)$.  Using this 
generator of $\pi_4(S^4)$, we orient the fiber $S^{3}$ with  
$\iota_3 \in  \pi_3(S^3)$ such that $\partial(\iota_3) = n \, \iota_4$ and
$n > 0$.  With this convention $M_{m,n}$ is an oriented manifold.   From
now on
all manifolds shall be smooth, compact and oriented and maps shall be
orientation preserving unless explicitly noted.  Thus ``$M'$ and $M$ are
homeomorphic'' shall
usually mean that $M'$ and $M$ are orientation preserving homeomorphic and
orientation
reversing maps will be considered explicitly as maps from $M'$ to $-M$.

\subsection{Structure of proof}

We wish to emphasize that our classification does not first solve the
problem of homotopy equivalence and then proceed to the finer
relations of homeomorphism and diffeomorphism.  Instead, we start 
by completing the $PL$-homeomorphism classification of
Wilkens.  We are then able to descend to give a homotopy classification by
explicitly exhibiting the entire normal invariant set of any $M_{m,n}$
as fiber homotopy equivalences.
\begin{Proposition}
\label{lem:structureset}
For all $j \in \Z$ there exist fiber homotopy equivalences
$f_j:M_{m+12j,n} \rightarrow M_{m,n}$ with normal invariant $j \in
\Z_{n}.$
\end{Proposition}
It follows from a simple application of surgery theory that a given
seven dimensional manifold $M'$ is homotopy  equivalent to $M_{m,n}$ if
and only if $M$ is homeomorphic to $M_{m+12j,n}$ for some $j$.  Given
the homeomorphism classification of the bundles this is enough to
yield the homotopy classification.

\section{Homeomorphism classification}

The elementary algebraic topology of our manifolds is determined by
the Euler number of the bundle $e(\xi_{m,n}) = n$.  From now on we 
use $H^{i}(M)$ to denote the $i$-th cohomology group of $M$ with 
integer coefficients.
$$\begin{array}{l}
H^{0}(M_{m,n}) \cong H^{7}(M_{m,n}) \cong \mathbf Z \\
H^{4}(M_{m,n}) \cong {\mathbf Z}_{n} \\
H^{i}(M_{m,n}) \cong 0 \ \mbox{ for all } \ i \ne 0, 4, 7. 
\end{array}$$

One sees from the Serre spectral sequence of the fibration $S^3
\hookrightarrow M_{m,n} \stackrel{\pi}{\rightarrow} S^4$ that
$\pi^*:H^4(S^4) \rightarrow H^4(M_{m,n})$ is a surjection.  Let
$\kappa_4 := \pi^*(\iota_4)$ and use $\kappa_4$ to identify $H^4(M_{m,n})$
with $\Z_n$.  The bundle structure ensures
that the linking form of $M_{m,n}\,,\,lk(M_{m,n}),$ is isomorphic to the
standard form, $l$, for all $m \in \Z$.  Thus,
$$\begin{array}{ccccl}
lk(M_{m,n}): & H^4(M_{m,n}) \cong \Z_{n}  \times  H^4(M_{m,n}) \cong \Z_{n} &
\longrightarrow &
\Q/\Z
\\ &(r  ,  s) & \longmapsto & l(r,s)=\frac{r\,s}{n} \,.
\end{array}$$

To see this, note that the linking form $lk(M_{m,n})$ is induced from the
intersection form of any coboundary $W$.  For
sphere
bundles we may choose $W=W_{m,n}$ to be the associated disk bundle,  
$D^4 \hookrightarrow W_{m,n} \stackrel{\pi_{W}}{\rightarrow}
S^{4}$,
$\p W_{m,n}= M_{m,n}$.  The long exact sequence of the pair
$(W_{m,n},M_{m,n})$
yields
$$0 \rightarrow H^4(W_{m,n},M_{m,n}) \stackrel{j^*}{\rightarrow}
H^4(W_{m,n}) \stackrel{i^*}{\rightarrow} H^4(M_{m,n}) \rightarrow 0.$$
Let $x = \pi_{W}^* \iota_4 \in H^4(W_{m,n}) \cong \Z$ and let
$y$ be a generator of $H^4(W_{m,n},M_{m,n}) \cong \Z$ such that $j^{*}(y) = nx$.
Since $\kappa_4$ = $i^*(x)$, the self linking number of the class $\kappa_4$ is
defined
as follows, 
$$lk(\kappa_4,\kappa_4) =\frac{1}{n}\,\<y \cup x, [W_{m,n},M_{m,n}]\> =
1.$$ 
Bilinearity is now enough to ensure that $lk(M_{m,n}) = l$.  Note that the
linking form $lk(M_{m,n})$ may be computed in a variety of ways.  See for
example \cite{KiSh} for a different method. 

\bigskip

  We now describe
the invariants required for our classification in some detail.

\subsection{The invariants $\mu$ and $s_{1}$}
The invariant $\,\mu(M_{m,n}) \equiv \mu(W,
M_{m,n}) ~ \mod 1\,$ is computed for any spin coboundary $W$
which satisfies the $\mu$-condition.  We choose
$W=W_{m,n}$ as above and note that the pair $(W_{m,n},M_{m,n})$ satisfies the
$\mu$-condition.  This coboundary was already used in \cite{M} and
\cite{T} and we use it to obtain for the $\mu$-invariant of the total 
spaces $M_{m,n}$
\begin{equation}\label{Inv2}
\begin{array}{lcl}
\mu(M_{m,n}) & \equiv & \frac{1}{2^{7} \cdot 7} \ \{p_{1}^{2}(W_{m,n}) - 4
\ 
\sigma[W_{m,n}]\} ~ \mod 1 ,
\end{array} 
\end{equation}
where $\sigma[W_{m,n}]$ is the signature of the quadratic form given by 
$$
\begin{CD}
H^{4}(W_{m,n}, \partial W_{m,n})  @>>> \Q \\
v  @>>> \< v \cup v, [W_{m,n}, \partial W_{m,n}] \>.
\end{CD}
$$
and $(p_{1}(W_{m,n}))^{2}$ is the characteristic number
$$\< (j^{*})^{-1} p_{1}(W_{m,n})  \cup p_{1}(W_{m,n}),[W_{m,n}, \partial 
W_{m,n}] \> $$
where $j^{*}:H^{4}(W_{m,n}, \partial W_{m,n}; \Q) {\longrightarrow} 
H^{4}(W_{m,n}; \Q)\,.$  As $n > 0$ we obtain for the signature
$\sigma(W_{m,n})=1$.
Furthermore, we chose $y$ such that
 $j^{*}(y) = n \,x$.  Hence $(j^{*})^{-1}(x)
= \frac{1}{n} \,y$.  Finally, it is well known (see \cite{M}) that
$\,p_{1}(W_{m,n}) = 2\,(n + 2m)\,x$ and thus 
\begin{equation}
\begin{array}{lcl}
 & &\<(j^{*})^{-1} p_{1}(W_{m,n})  \cup p_{1}(W_{m,n}),[W_{m,n}, \partial 
W_{m,n}] \> \\
& \rule{0cm}{0.7cm} = & \< \frac{1}{n}\,2\,(n + 2m)\,x \cup 2\,(n + 
2m)\,x, [W_{m,n}, \partial W_{m,n}] \> \\
& \rule{0cm}{0.7cm} = & \frac{4\,(n + 2\,m)^{2}}{n}\,.
\end{array} 
\end{equation}
Hence 
$$\mu(M_{m,n}) \equiv \frac{1}{2^{5} \cdot 7 \cdot n}\,\{(n + 
2\,m)^{2} - n\} ~ \mod 1\,$$
and therefore 
\begin{equation}
\begin{array}{lcl}
s_{1}(M_{m,n}) = 28\mu(M_{m,n}) & \equiv &\frac{1}{8 \cdot n}\,\{(n + 
2\,m)^{2} - n\} ~ \mod 1 \\
& \rule{0cm}{0.7cm}\equiv & \frac{m^{2}}{2n} + 
\frac{m}{2} + \frac{n-1}{8}  ~ \mod 1 \,.
\end{array} 
\end{equation}
\noindent
As noted above, it is proven in \cite{KS} that $s_1$ is an invariant of
topological spin manifolds.

\medskip

We would like to point out that there are mistakes in the literature in 
calculations of the $\mu$-invariant for the manifolds $M_{m,n}$.  
In particular, an incorrect formula for the $\mu$-invariant \cite[p.
67]{D}
is applied to prove \cite[Theorem 2.1]{D}.  In fact, the statement 
\cite[Theorem 2.1]{D} is true under use of the true $\mu$-invariant as 
W. Ziller pointed out to us.  However, Ziller tells us that the proof 
of \cite[Theorem A]{Ma} does not work with the true $\mu$-invariant.  
We do not know if \cite[Theorem A]{Ma} is true.

\subsection{The invariant $\halfp$}
Since the total spaces $M_{m,n}$ are $2$-connected they all have unique
spin structures and we work with spin characteristic classes for 
simplicity.  The Hopf bundle $\pi_{0,1}: \xi_{0,1} \longrightarrow 
S^{4}$ defines a generator of $\pi_{4}(BSpin)$ but $p_{1}(\xi_{0,1}) = 
\pm 2 \in H^{4}(S^{4})$.  It follows that
half the first Pontrjagin class, $\halfp$, is a generator of
$H^4(BSpin)$ and thus is always integral.   Moreover,
$\halfp$ is a topological invariant since the canonical map $H^4(BSpin)
\rightarrow H^4(BTop)$ is an injection.  We refer the reader to 
\cite[Lemma  6.5]{KS} for a discussion of this map.  In \cite{M}
Milnor showed how to compute $p_{1}$ for the manifolds $M_{m,n}$.  One
first observes that the stable tangent bundle of a sphere bundle is a
pullback of bundles over the base space.  As we shall need the
corresponding
fact for the stable normal bundle we give both results now.  Let
$\tau_{B}$
denote  the stable tangent bundle of the base space $B$ and let $\nu_{B}$ 
denote the corresponding stable normal bundle.  
\begin{fact} Let $S(\psi)$ be the sphere bundle associated to a vector
bundle
$\psi$ over a manifold $B$ and let $-\psi$ denote the stable inverse of
$\psi$.  Then
$$\nu_{S(\psi)} = \pi_{\psi}^*(\nu_B \oplus -\psi)
~\mbox{and}~~\tau_{S(\psi)} = \pi_{\psi}^*(\tau_B \oplus \psi).$$
\label{fact}
\end{fact}
Recalling that $p_{1}(\xi_{m,n})= 2\,(n+2m)\,\iota_4$ we obtain that
\begin{equation}
\begin{array}{lcl}
\halfp(M_{m,n}) & = & \halfp (\tau_{M_{m,n}})
= \halfp (\pi^*(\tau_{S^4} \oplus \xi_{m,n}))
= \pi^*(\halfp(\xi_{m,n}))
= \pi^*((n+2m)\,\iota_4) \\
& \rule{0cm}{0.7cm} = & 2\,m \,\kappa_4 \in \Z_n \,.
\end{array}
\end{equation}

The characteristic class $\halfp$ is identical to the invariant
$\hat\beta$
defined in \cite{Wi} as the  obstruction to the tangent bundle being
trivial
over the $4$-skeleton.  Wilkens' full invariant for a $2$-connected
$7$-manifold with torsion fourth cohomology group is the triple
$(H^4(M),lk(M),\halfp(M))$ where $lk(M)$ denotes the linking form of $M$
on
$H^4(M)$.  He denoted such a triple abstractly as
$(G,l,\beta)$ and defined two triples $(G,l,\beta)$ and $(G',l',\beta')$
to be 
the same invariant if 
there is an isomorphism $h:G \rightarrow G'$ preserving the linking form
and sending $\beta$ to $\beta'$.  For $M=M_{m,n}$ we have identified
$H^4(M_{m,n})$ with $\Z_n$ by taking $\pi_{m,n}^* \iota_4  = \kappa_4$ as
a
generator of $H^4(M_{m,n})$.  With this identification the linking form 
becomes $lk(M_{m,n}) = l$ where $l(r,s) = \frac{rs}{n}$.  Hence when we
speak of
the invariant $\halfp$ we interpret it as part of the Wilkens triple
$(\Z_n,l,\halfp(M_{m,n})=2m)$.  Note that the automorphisms of 
$\Z_{n}$ which preserve the linking form $l$ can be described as 
elements of the set $A^{+}(n):= \{\alpha \in \Z_{n}: \alpha^{2} = 1\}.$   
Therefore we say that $M_{m',n}$ and $M_{m,n}$ have the same $\halfp$
when there is an $\alpha \in A^+(n)$ such that
$\halfp(M_{m',n}) = \alpha \,\halfp(M_{m,n}) \in \Z_n$.  In the 
orientation reversing case the invariant $\halfp$ is unchanged but for 
the linking form we obtain 
$lk(-M_{m,n}) =
-lk(M_{m,n}) = -l$.  We must therefore consider the automorphisms of
$\Z_{n}$
which send $l$ to
$-l$.  These form the set $A^{-}(n):= \{\alpha \in \Z_{n}: \alpha^{2} = -
1\}$.
Thus we say that two manifolds $M_{m',n}$ and $-M_{m,n}$ have 
the same $\halfp$ when there exists an $\alpha \in A^-(n)$ such that
$\halfp(M_{m',n}) = \alpha \, \halfp(M_{m,n}) \in \Z_n$.
Elementary considerations in  number theory lead to the following lemma 
describing the automorphisms $\alpha$, see for
example 
\cite[Theorem 5.2]{L}.

\begin{Lemma}\label{numbertheory} 

\begin{description}
	
\item{{\bf (a)}}  The congruence $\alpha^{2} \equiv 1 ~ \mod n$ has 
$2^{r+u}$ solutions, where $r$ is the number of distinct odd prime 
divisors of $n$ and $u$ is $0, 1$ or $2$ according to $4$ does not 
divide $n$, 
$2^{2}$ exactly divides $n$ or $8$ divides $n$.  
\item{{\bf (b)}}  The congruence $\alpha^{2} \equiv - 1 ~ \mod n$ is 
solvable if and only if 
$n=p_1^{i_1} \dots p_k^{i_k},~p_i \equiv 1 ~ \mod 4\,, p_{i}$ 
prime or $n=2p_1^{i_1} \dots p_k^{i_k},~p_i \equiv 1 ~ \mod 4, p_{i}$ 
prime.  In this case one has $2^{k}$ solutions.
\end{description}
\end{Lemma}

\subsection{Orientation Preserving Homeomorphism}
\label{subsec:top}
Let $G^{*}$ be the torsion subgroup of $G$.  Then Wilkens called a 
symmetric bilinear form $b:G^{*} \times G^{*} \rightarrow
\Q/\Z$ \emph{irreducible} if it cannot be written as the proper sum of two
bilinear forms and he defined $M$ to be \emph{indecomposable} if $H^4(M)$
is
finite and $lk(M)$ is irreducible or if $H^4(M) \cong \Z$.  Wall
(\cite{Wa2}) had already proven that finite irreducible forms only exist
for
$G \cong \Z_{p^a}$ or $\Z_{2^a} \oplus
\Z_{2^a}$ for $p$ prime.  Modulo a $\Z_2$ ambiguity in some cases Wilkens
showed 
that the triple $(H^4(M),lk(M),\halfp(M))$ classifies indecomposable
$2$-connected $7$-manifolds up to oriented diffeomorphism
and  connected sum with an exotic $7$-sphere.  In the case of
$M=M_{m,n}$ we obtain
$H^4(M_{m,n}) \cong \Z_{n} \cong  \bigoplus_{i=1}^{n} \Z_{q_{i}^{e_{i}}}$
where $n = \prod_{i=1}^{k} q_{i}^{e_{i}}$ is the prime decomposition 
of $n$.  Then any
automorphism of $H^4(M_{m,n})$ must preserve the decomposition into
irreducible
forms, there being no nontrivial homomorphisms between cyclic groups of
coprime
orders.  Hence any homeomorphism $h:M_{m',n}\rightarrow M_{m,n}$ must
preserve the
decompositions of $M_{m',n}$ and $M_{m,n}$ into indecomposable manifolds.
It
follows that $\halfp(M_{m,n}) = (\Z_n,l,2\,m)$ classifies the total spaces
$M_{m,n}$ up to an occasional $\Z_2$ ambiguity. 

A homeomorphism of manifolds $f:M' \rightarrow M$ is called an almost
diffeomorphism if there is an exotic sphere $\Sigma$ such that $f:M'
\rightarrow M \,\sharp \,\Sigma$ is a diffeomorphism.  As $PL/O$ is
$6$-connected, general smoothing theory implies that the manifolds
$M_{m,n}$ are
almost-diffeomorphic if and only if they are $PL$-homeomorphic \cite{MM}.
Below
we show that Wilkens' ambiguous cases have distinct values for $s_1$.  It
follows
that manifolds $M_{m',n}$ and $M_{m,n}$ are $PL$-homeomorphic if and only
if they
have the same invariants $\halfp$ and $s_1$.  But these are topological
invariants and so $M_{m',n}$ and $M_{m,n}$ are homeomorphic if and only if
they
are $PL$-homeomorphic.

\begin{Remark}
When $n$ is even, $H^4(M_{m,n};\Z_{2}) \cong \Z_2$ and thus
Kirby-Siebenmann theory guarantees that there is a pair of non-concordant 
$PL$-structures on $M_{m,n}$.  The above argument shows that these $PL$
structures are $PL$-homeomorphic. 
\end{Remark}

Let $n =2^{a}\,q$ with  $q$ odd.  It follows directly from Wilkens that
$\halfp$
classifies the manifolds $M_{m,n}$ in the following circumstances: $n$ odd
or
$n=4\,q~$ and $m$ odd or $n = 2^{a}\,q\,,\,a>2$ and $m$ even.  However,
there
remained a $\Z_2$ ambiguity in the cases: $n = 2\,q$ with any 
$m$, $n = 4\,q~$ with  $m$ even and
$a>2\,$ with $m$ odd.  We now complete the orientation preserving
homeomorphism classification by  showing that the invariant $s_{1}$
distinguishes
manifolds $M_{m,n}$ with the same $\halfp$ invariant in the remaining
cases.  Let
$m' = m - \frac{n}{2}$. Then $\halfp(M_{m',n}) = \halfp(M_{m,n})$ and
\begin{equation}
\begin{array}{lcl}
 s_{1}(M_{m,n}) - s_{1}(M_{m',n})& \equiv & \frac{m - m'}{2} + 
 \frac{m^{2} - m'^{2}}{2n} ~ \mod 1 \\
 & \rule{0cm}{0.7cm}\equiv & \frac{m}{2} + \frac{n}{8}  ~ \mod 1 \\
  & \rule{0cm}{0.7cm}\equiv & \frac{m}{2} + 
  \frac{2^{a}q}{8} ~ \mod 1 \not\equiv 0 ~ \mod 1
\end{array} 
\end{equation}
for the cases  of $a=1$ or $a=2$ and $m$ even or $a>2$ and $m$ odd.  
We have thus proven Theorem \ref{top} in the orientation preserving 
case.   Now, Wilkens showed that every $\halfp(M)$ of a $2$-connected
$7$-manifold must be divisible by two and since $\halfp(M_{m,n})=2\,m$,
$S^3$-bundles over $S^4$ realize every Wilkens triple
$(H^4(M),lk(M),\halfp(M))$
with linking form $lk(M) = l, l(r, s) = \frac{rs}{n}\,.$  
Moreover, we have just seen that each of Wilkens' ambiguous cases
is realized by such a bundle.  It follows then from Wilkens'
classification that any $2$-connected $7$-manifold with linking
form  $l$ is almost diffeomorphic, and
hence homeomorphic, to an $S^3$-bundle  over $S^4$.  This is
precisely the statement of Corollary \ref{topcor2}.  We now proceed to
demonstrate the algebraic formulation of Theorem \ref{top} given in
Corollary
\ref{topcor}.

The first two cases of the orientation preserving homeomorphism 
condition of Corollary \ref{topcor} immediately follow from 
the fact that in those cases the invariant $\halfp$ is enough to 
classify the manifolds.  For the last case we have to prove that the 
conditions
$\halfp(M_{m',n}) = \alpha \, \halfp(M_{m,n})$ and $s_{1}(M_{m',n}) =
s_{1}(M_{m,n})$ where $\alpha^2 \equiv \,1 ~ \mod n$ are equivalent 
to the conditions
$m' \equiv \alpha \, m ~\mod n$ and $\alpha \equiv \pm \, 1 ~\mod 2^a$ 
where  $\alpha^2 \equiv \,1\, ~ \mod n$.  Let $\alpha^{2} = k\,n +1$ 
for some $k \in \Z$ and let $m' \equiv \alpha \,m + \epsilon \, 
\frac{n}{2} ~\mod n $ where $\epsilon = 0, 1$.  Then we calculate the 
difference of the $s_{1}$ invariants:
$$\begin{array}{lcl}
s_{1}(M_{m,n}) - s_{1}(M_{m',n}) &\equiv& \frac{\alpha m - m}{2} + 
\frac{\epsilon n}{4} + \frac{(\alpha m + \frac{\epsilon n}{2})^{2}}{2n} - 
\frac{m^{2}}{2n} ~ \mod 1\\
&\rule{0cm}{0.7cm}\equiv& \frac{k m^{2}}{2} + \frac{3 \epsilon n}{8} +
\frac{\alpha m 
\epsilon}{2} ~\mod 1\,.
\end{array} $$

\noindent
Now if $m' \equiv \alpha \, m ~\mod n$ and $\alpha \equiv \pm \, 1 
~\mod 2^a$ 
where  $\alpha^2 \equiv \,1\, ~ \mod n$, then the $\halfp$ 
invariants must be the same.  But the $s_{1}$ invariants also must
coincide as $\alpha \equiv \pm \, 1 ~\mod 2^a$ implies that $k$ 
is even and the $\halfp$ condition implies that $\epsilon = 0$.  For 
the converse we argue case by 
case.  In the first case of $a=1$ we need to show that 
$m' \equiv \alpha \, m ~\mod 2q$ where $\alpha \equiv 1 ~\mod 2$.  But
here 
$0 \equiv s_{1}(M_{m,n}) - s_{1}(M_{m',n}) \equiv 
\frac{k m^{2}}{2} + \frac{3 \epsilon q}{4} + \frac{\alpha m 
\epsilon}{2} ~\mod 1$ which implies that $\epsilon = 0$ and hence 
$m' \equiv \alpha \, m ~\mod 2 q$.  Also $\alpha^2 \equiv \,1\, ~ \mod 
2q$ reduces modulo $2$ to 
$\alpha^2 \equiv \,1\, ~ \mod 2$ and hence $\alpha \equiv 1 ~\mod 2$. 
In the second case of $a=2\,,\,m$ 
even we need to show that  $m' \equiv \alpha \, m ~\mod 4q$ and 
$\alpha \equiv \,\pm 1 ~\mod 4$.  Again  
$0  \equiv 
\frac{k m^{2}}{2} + \frac{3 \epsilon q}{2} + \frac{\alpha m 
\epsilon}{2} \equiv \frac{3 \epsilon q}{2} ~\mod 1  $ implies that 
$\epsilon = 0$ and hence 
$m' \equiv \alpha \, m ~\mod 4 q$.   Also $\alpha^2 \equiv \,1\, ~ \mod 
4q$ reduces modulo $4$ to 
$\alpha^2 \equiv \,1\, ~ \mod 4$ and hence $\alpha \equiv \pm 1 ~\mod 4$. 
 For the last case of $a > 2\,,\,m$ 
odd we need to show that  $m' \equiv \pm \, \alpha ~\mod 2^{a}q$ 
and $\alpha \equiv 1 ~\mod 2^{a}$.  Again
$0  \equiv 
\frac{k m^{2}}{2} + \frac{3 \epsilon 2^{a} q}{8} + \frac{\alpha m 
\epsilon}{2} \equiv \frac{k m^{2}}{2}  + \frac{\alpha m 
\epsilon}{2} \equiv \frac{k + \epsilon}{2} ~\mod 1  $ as $m$ is odd.  
By Lemma \ref{numbertheory} we know that there are $4$  
possible cases for $\alpha$ after reduction modulo $2^{a}$, namely 
$\alpha \equiv 1, 2^{a-1} - 1, 2^{a-1} + 1,
2^{a} - 1 ~ \mod 2^{a}\,,\,a >2\,.$  
For $\alpha \equiv 2^{a-1} - 1, 2^{a-1} + 1 ~ \mod 2^{a}$ we obtain 
that $k$ must 
be odd as $\alpha^{2} \equiv 1 ~\mod n$.  Hence $\epsilon$ must be 
odd as well by the $s_{1}$ condition.  But then we obtain that 
$m' \equiv \pm \, m ~\mod 2^{a}$ and hence 
$m' \equiv \bar\alpha \, m ~\mod n$ with 
$\bar\alpha \equiv \pm \, 1 ~\mod 2^{a}\,.$  
In the case 
$\alpha \equiv 1,2^{a} - 1 ~ \mod 2^{a}$ we obtain that $k$ must be 
even and hence $\epsilon = 0$ by the $s_{1}$ condition.  Again
we conclude that $m' \equiv \alpha  \, m ~\mod n$ and 
$\alpha^{2} \equiv \pm 1 ~\mod 2^{a}\,.$

\subsection{Orientation Reversing Homeomorphism}

There is an orientation reversing homeomorphism between 
manifolds $M$ and $M'$ exactly when there is an orientation 
preserving homeomorphism between $M$ and $-M'$.  Since the linking form of
$-M_{m,n}$ is $-l$  we now seek automorphisms of $\Z_{n}$ which 
send $l$ to $-l$.  Such automorphisms $\alpha$
are precisely the elements of $A^-(n) \subset \Z_n$.   Now $\halfp$ is
unchanged under changes of orientation whereas $s_1(-M) = -s_1(M)$. 
We conclude  that $M_{m',n}$ is orientation reversing homeomorphic to
$M_{m,n}$ if and only if $s_1(M_{m',n}) = -s_1(M_{m,n})$ and
$\halfp(M_{m',n}) = \alpha \, \halfp(M_{m,n})$ for some $\alpha \in A^-(n)$.
Applying Lemma \ref{numbertheory}
 shows that $A^-(n)$ is nonempty only if
$n=2^a \, q$ where $a = 0$ or $1$ and $q$ is a product of powers of primes
which
are congruent to $1~\mod 4$.  When $a=0$, $\halfp$ alone
classifies the manifolds $M_{m,n}$ and Theorem \ref{top} follows in the
orientation reversing case.  When 
$a=1$, $\halfp$ is ambiguous and we resort to $s_1$ to settle the
issue.  Assume that $\halfp(M_{m',n})=\alpha \, \halfp(M_{m,n})$ so that $m' 
\equiv 
\alpha m + \epsilon q ~ \mod n$ where $\epsilon = 0$ or $1$.  One
calculates
$$s_1(M_{m',n}) + s_1(M_{m,n}) = \frac{\epsilon}{2} +
\frac{\epsilon^2}{4} +\frac{1}{4} ~ \mod 1. $$
Thus $s_1(M_{m',n}) \equiv - s_1(M_{m,n}) ~ \mod 1$ if and only if
$\epsilon =
1$ and the orientation reversing case of Theorem \ref{top} follows.

\section{Diffeomorphism classification}

From section  \ref{subsec:top} we know that the manifolds $M_{m,n}$ and 
$M_{m',n}$ are 
$PL$-homeomorphic if and only if they are almost diffeomorphic.  
In  \cite{EK} it is shown that the 
invariant $\mu$ is additive with respect to connected sums and that it
distinguishes all exotic $7$-spheres.  Combining these 
facts with the $PL$-homeomorphism
classification we obtain that the total spaces $M_{m,n}$ and $M_{m',n}$
are 
diffeomorphic if and only if 
they are $PL$-homeomorphic and their $\mu$-invariants coincide.  
Hence they are diffeomorphic if and only if they have the same 
invariants $\halfp$ and $\mu$.

Smooth surgery theory (see \cite{MM}) implies that there are exactly 
$28$ different smooth manifolds homeomorphic to $M_{m,n}$.  A natural 
question to ask is which manifolds homeomorphic to $M_{m,n}$ may be 
represented by total spaces of $S^{3}$-bundles over $S^{4}$, i.e. by 
manifolds $M_{m',n}$.  In order to answer this question for some 
examples of small $n$ we introduce the following integer valued 
functions.

\begin{Definition}\label{intfct}
\begin{enumerate}
\item  Let $\h^{+}(n) ~[\h(n)]$ be the number of orientation 
preserving [orientation preserving and reversing] homeomorphism types 
of the total space of an $S^{3}$-bundle over $S^{4}$ with Euler class $n$.
\item  Let $\Diff^{+}(m,n) ~[ \Diff(m,n)]$ be the number of distinct 
orientation preserving [orientation preserving and reversing] classes 
of smooth manifolds represented by $S^{3}$-bundles over $S^{4}$ in 
each homeomorphism class of $M_{m,n}$.
\end{enumerate}
\end{Definition}

Using the congruences of Theorems \ref{hom} 
and \ref{top} we obtain the following results for a selection of 
values of $n$ with $n \le 16$.  Note that in this range orientation 
reversing homeomorphism and diffeomorphisms only exist for $n = 1, 2, 
5, 10, 13\,.$  

\begin{itemize}

\item $\h^{+}(1) = 1$  and $\Diff^{+}(m,1) = 16\,;$

      $\h(1) = 1$  and $\Diff(m,1) = 11\,;$
	  
\smallskip 

\item $\h^{+}(2) = 2$  and $\Diff^{+}(m,2) = 8$ in each homeomorphism 
class;

$\h(2) = 1$  and $\Diff(m,2) = 13\,$ in each homeomorphism class;

\smallskip 

\item $\h^{+}(5) = 3$  and $\Diff^{+}(m,5) = 16$ in each 
homeomorphism class; 

$
\begin{array}{lll} 
 \h(5) = 2 ~ \mbox{ and } & \Diff(m,5) & =  12 ~ \mbox{for} ~ m \equiv ~ 0 ~ \mod 
 5\,;\\
& \Diff(m,5) & =  24 ~ \mbox{for} ~ m \not\equiv ~ 0 ~ \mod 5\,;  
\end{array}
$

\smallskip 

\item $\h^{+}(10) = 6$  and $\Diff^{+}(m,10) = 8$ in each homeomorphism 
class;

      $\h(10) = 3$  and $\Diff(m,10) = 14\,$ in each homeomorphism 
class;
 
\smallskip 

\item $
\begin{array}{lll} 
 \h^{+}(7) = 4 ~ \mbox{ and } & \Diff^{+}(m,7)  =  4 & ~ \mbox{for} ~ m \equiv ~ 0 ~ \mod 
 7\,;\\
& \Diff^{+}(m,7)  =  28 & ~ \mbox{for} ~ m \not\equiv ~ 0 ~ \mod 7\,;  
\end{array}
$

\smallskip 

\item $
\begin{array}{lll} 
 \h^{+}(14) = 8 \, \mbox{ and } & \Diff^{+}(m,14)  =  2 & ~ \mbox{for} \, m \equiv ~ 0 ~ \mod 
 14 ~\mbox{and} \\
 &  & ~ \mbox{for} \, m \equiv \, 7 \, \mod 
 14\, ;\\
& \Diff^{+}(m,14)  =  14 & ~ \mbox{for all others;}   
\end{array}
$

\smallskip 

 \item $
\begin{array}{lll} 
 \h^{+}(4) = 3 \, \mbox{ and } & \Diff^{+}(m,4)  =  4 & ~ \mbox{for} \, m \equiv ~ 0 ~ \mod 
 4 \, ;\\
& \Diff^{+}(m,4)  =  8 & ~ \mbox{for}   \, m \equiv ~ 2 ~ \mod 
 4 \, ;\\
 & \Diff^{+}(m,4)  =  16 & ~ \mbox{for}   \, m \equiv ~ \pm 1 ~ \mod 
 4 \, ;
 \end{array}$
 
 \smallskip 
 
  \item $
\begin{array}{lll} 
  \h^{+}(8) = 4 \, \mbox{ and } & \Diff^{+}(m,8)  =  8  & ~ \mbox{for} \, 
  m \equiv ~ 0 \, \mbox{ or } \, 4 ~ \mod 8 \, \mbox{ and }\\
&   & ~ \mbox{for}   \, m \equiv ~ 2 ~ \mod 
 8 \, ;\\
 & \Diff^{+}(m,8)  =  16 & ~ \mbox{for}   \, m \equiv ~ \pm 1, \pm 3 ~ \mod 
 8 \, ;
 \end{array}$
 
 \smallskip 
 
 \item $\begin{array}{lll} 
 \h^{+}(12) = 6 \, \mbox{ and }  & \Diff^{+}(m,12)  =   4 & ~ \mbox{for} \, 
     m  \equiv  ~ 0, \pm 4 ~ \mod 12 \,;\\
 &\Diff^{+}(m,12)  =  8 & ~ \mbox{for}   \, m  \equiv  ~ \pm 2 , 6 ~ \mod 
 12 \,; \\
  & \Diff^{+}(m,12)   =  16 & ~ \mbox{for}   \, m  \equiv  ~ \pm 
  1\,\mbox{or} \, \pm 5 ~ \mod 
 12 \,\mbox{ and }\\
 &  &  ~ \mbox{for} ~  m  \equiv  ~ \pm 3 ~ \mod 12  ;
 \end{array}$
 
 \smallskip 
 
 \item $\begin{array}{lll} 
 \h^{+}(16) = 7 \, \mbox{ and }  &\Diff^{+}(m,16)  =   4  & ~ \mbox{for} \, 
  m  \equiv  ~ \pm 4 ~ \mod 16 \,;\\
 &\Diff^{+}(m,16)  = 8 & ~ \mbox{for}   \, m  \equiv  ~ 0 \, 
 \mbox{or} \, 8 ~ \mod 16 \, ;\\
 &\Diff^{+}(m,16)  =  16 & ~ \mbox{for}  \, m  \equiv  ~ \pm 
  2\,\mbox{or} \, \pm 6 ~ \mod 
 16 \,\mbox{ and }\\
 &  &  ~ \mbox{for} ~  m  \equiv  ~ \pm1, \pm 3, \pm 5, \pm 7 ~ \mod 
 16\,.\\
 
\end{array}
$
\end{itemize}

\section{Classification up to homotopy equivalence}

In order to use surgery
theory to classify these manifolds we recall the
following definitions based on the work of Sullivan, see \cite{MM} for
a detailed description.  Let $M$ be a $\Cat$-manifold, where 
$\Cat$ stands for the piecewise linear ($\PL$), topological
($\Top$), or smooth  (${\mathrm O}$) category.  Then ${\mathcal 
S}^{\Cat}(M)$,
the $\Cat$-structure set of $M$, consists of pairs $(L, f)$ where
$L$ is a $\Cat$-manifold and $f: L \rightarrow M$ is a simple
homotopy equivalence.  Two objects $(L_{1}, f_{1})$ and $(L_{2},
f_{2})$ represent the same element in ${\mathcal S}^{\Cat}(M)$ if
there exist an h-cobordism $W$ with $\partial W = L_{1} \cup L_{2}$
and a homotopy equivalence $F: W \rightarrow M$.  Moreover we denote
by ${\mathcal N}^{\Cat}(M)$ the set of equivalence classes of normal
cobordisms.  Here two normal maps (or surgery problems)
$$
\begin{CD}
\nu_{L_{1}}  @>\hat{f_{1}}>> \nu_{1} @. \phantom{xxxxxxxxxxxxxxx}  
\nu_{L_{2}} @>\hat{f_{2}}>> \nu_{2}\\
@VVV          @VVV   \phantom{xxxxxx} \mbox{ and }  \phantom{xxxxxx}@VVV
@VVV    \\
L_{1} @>f_{1}>> M                @.  \phantom{xxxxxxxxxxxxxxx}  
L_{2} @>f_{2}>> M
\end{CD}
$$
are normally cobordant if there exists a cobordism  $W$ with 
$\partial W = L_{1} \cup L_{2}$ and 
$$
\begin{CD}
\nu_{W}  @>\hat{F}>> \nu \times I \\
@VVV          @VVV         \\
W @>F>> M \times I    
\end{CD}
$$
such that $F|L_{1} = f_{1}\,,\,\hat{F}|\nu_{L_{1}} = \hat{f_{1}},\,
F|L_{2} = f_{2}\,,\,\hat{F}|\nu_{L_{2}} = b \circ \hat{f_{2}}$ where 
$b: \nu_{2} \rightarrow \nu_{1}$ is a bundle isomorphism.  
Let 
$$ L_{n} = \left.\{ \begin{array}{ll}
	   {\mathbf Z}, & \mbox{ if } n = 4\,k;\\
	   {\mathbf Z}_{2} , & \mbox{ if } n = 4\,k + 2;\\
	    0 , & \mbox{otherwise.} 
	\end{array} 
	\right.\} $$
The following two facts are well-known, see  \cite{MM}.
\medskip

\noindent 
(1)  ${\mathcal N}^{Cat}(M) \cong [M, G/\Cat]\,.$
\medskip

\noindent
(2)  {\it There is an exact sequence of sets, the surgery exact sequence:}
$$
\cdots \longrightarrow {\mathcal N}^{\Cat}(M^{n} \times I, \partial) 
\stackrel{s}{\longrightarrow}
L_{n+1} \longrightarrow {\mathcal S}^{\Cat}(M^{n}) \longrightarrow
{\mathcal N}^{\Cat}(M^{n}) \longrightarrow L_{n}.
$$

A natural set of homotopy equivalences between sphere bundles is the set
of fiber
homotopy equivalences of the associated spherical fibrations.  The
transition from
vector bundles to spherical fibrations is encoded in the $J$-homomorphism,
$J_{n,k}:\pi_n(SO(k))
\rightarrow \pi_{n+k}(S^k)$,
see \cite[pp 502-4]{Wh} for a standard definition.
If $\zeta$ and $\xi$ are rank $k$ vector bundles over $S^n$, $\zeta, \xi
\in
\pi_{n-1}(SO(k))$, then $S(\zeta \oplus \epsilon)$ is fiber homotopically
equivalent to $S(\xi \oplus  \epsilon)$ if and only if $J_{n-1,k}(\zeta) =
J_{n-1,k}(\xi)$.  Here $S(\xi)$ denotes the sphere bundle of $\xi$
and $\epsilon$ is a trivial bundle.  In fact, it
is easy to see that one has fiber homotopy equivalence even without 
stabilization.

\begin{Lemma}
\label{lem:hes}
There exist fiber homotopy equivalences $f_j:M_{m+12j,n} \rightarrow
M_{m,n}$ for all $ j \in \Z$ .
\end{Lemma}

\noindent 
{\bf Proof.}
For our manifolds $n = 4$ and $k = 4$ and the $J$-homomorphism becomes 
$J_{3,4}:\Z \oplus  \Z \cong \pi_3(SO(4))
\rightarrow \pi_7(S^4) \cong \Z \oplus  \Z_{12}$.  Moreover the projection
of
$J_{3,4}$ onto the free summand is given by taking the Euler number of
the corresponding bundle. 
Hence $J_{3,4}(\xi_{m+12j,n})=J_{3,4}(\xi_{m,n})$ from which it follows
that  
$S(\xi_{m+12j,n})=M_{m+12j,n}$ and $S(\xi_{m,n})=M_{m,n}$ are fiber
homotopy equivalent. $\qed$ 

\vspace{0.3cm}

\noindent
We now give a general result which allows one to identify the normal
invariant defined by any fiber homotopy equivalence of sphere bundles.
\begin{Lemma}
\label{normalinvt}
Let $\zeta$ and $\xi$ be vector bundles over a simply connected closed
smooth manifold $B$ and let $f:S(\zeta) \rightarrow S(\xi)$ be a fiber
homotopy
equivalence of sphere  bundles over $B$.  Then $f$ defines an element
$[f]$ in
$\NO(B)$ and also,  as a homotopy equivalence, an element $\{f\}$ of
$\SO(S(\xi))$.
Let
$\pi_{\xi}^*:\NO (B) \rightarrow \NO (S(\xi))$ be the map induced by 
$\pi_{\xi}: S(\xi) \rightarrow B$ and let $\eta^{O}$ be the map from 
the structure set to the normal
invariant set.  Then,
$$\eta^{O}(\{f\}) = \pi_{\xi}^*([f]).$$
\end{Lemma}

\noindent 
{\bf Proof.}
The proof requires translating between alternate definitions of the
elements of the set of normal invariants of $S(\xi)$. 
See \cite[Chapter 1]{MM} for a detailed definition of the
normal invariant set.  
We begin with the homotopy equivalence $f$ which defines a class
$\{ f \}$ in the structure set of $S(\xi)$ and therefore a normal
invariant,
$\eta^O(\{ f\})$, of $S(\xi)$.  This normal invariant is the pullback 
of the stable
normal bundle of $S(\xi)$, $\nu_{S(\xi)}$, along a homotopy inverse of
$f$, denoted by $g$.  Identifying $f^*g^*(\nu_{S(\xi)})$ with
$\nu_{S(\xi)}$ yields the following description of $\eta^O( \{ f\})$.
$$\eta^O( \{ f\}) = \left( \begin{array}{crclcc}
\label{arr:henormalinvt}
& \nu_{S(\zeta)} & \stackrel{\cong}{\longrightarrow} &
f^*g^*(\nu_{S(\zeta)}) &
\longrightarrow & g^*\nu_{S(\xi)} \\
& \searrow & & \swarrow & & \downarrow \\
& & S(\zeta) & & \stackrel{f}{\longrightarrow} & S(\xi)
\end{array} \right) $$
If we choose $g$ to be a fiber homotopy equivalence itself then
$f^*g^*=id$
as all bundles
involved are pullbacks (recall Fact \ref{fact}).  Moreover,
$g^*\pi_{\zeta}^* =
\pi_{\xi}^*$.  We shall need explicit coordinates to compare
normal invariants.  Recall that for any pair of vector bundles $\phi$
and $\psi$,
$\pi_{\psi}^*(\phi) = \{ (v,w)| \pi_{\psi}(v) = \pi_{\phi}(w) \} \subset
S(\psi) \times \phi$.  With these coordinates,
\begin{equation}
\label{eq:nof}
\eta^O(\{ f\}) = \left( \begin{array}{lcclr}
& \pi_{\zeta}^*(\nu_B \oplus - \zeta) & \stackrel{g^*}{\longrightarrow} &
\pi_{\xi}^*(\nu_B \oplus - \zeta) & \\
& \downarrow & & \downarrow & ,\,g^{*}:(z,(y,z')) \mapsto (f(z),(y,z')) \\
& S(\zeta) & \longrightarrow & S(\xi) & \\
\end{array} \right) 
\end{equation}
\noindent
with $z \in S(\zeta)_b, y\in \nu_b, z' \in -\zeta_b$ and $b \in B$.  

Next we turn to the element $[f] \in \NO (B)$ which is the class of fiber
homotopy equivalences between $\zeta$ and $\xi$ containing the fiberwise
cone on $f$,
$f':\zeta \rightarrow \xi$.  The element $\pi_{\xi}^*([f])$ is the
pullback
of this
fiber homotopy equivalence class over $S(\xi)$.  
$$\begin{array}{ccccc}
& \pi_{\xi}^*(\zeta) & \stackrel{\pi_{\xi}^*f'}{\longrightarrow} & 
\pi_{\xi}^*(\xi) &  \\ 
& ~~~\searrow & & \swarrow ~~~  & ,\,\pi_{\xi}^*f':(x,z) \mapsto 
(x,f'(z)) \\ 
S(\zeta) & \stackrel{f}{\longrightarrow} & S(\xi) & \\
~\searrow & & \swarrow ~~~ & &\\
& B & & &
\end{array}$$

It remains to convert the fiber homotopy equivalence $\pi_{\xi}^*f'$ 
into a degree one normal map.  To do this we give a more general version
of the 
path from a fiber homotopy equivalence to a normal invariant than found in
\cite{MM}.

\begin{Lemma}
\label{obs:fhetod1nm}
Let $h:\phi \longrightarrow \psi$ be a fiber homotopy equivalence of
vector bundles over a manifold $M$, let $s:M \hookrightarrow \psi$ 
be a section of $\pi_{\psi}$ such that $h$ is transverse to $s(M)$ and let
$N:=h^{-1}(s(M))$ be the inverse image of $s(M)$.  Then the normal
invariant of $M$
associated to $h$ is 
$$\eta^O(h) = \left(\begin{array}{ccc}
\nu_N \cong (\pi_{\phi}|_N)^*(\psi \oplus \nu_M \oplus -\phi) &
\longrightarrow & \psi \oplus \nu_M \oplus \phi  \\
\downarrow & & \downarrow \\
N & \stackrel{\pi_{\psi}|N}{\longrightarrow} & M \\
(w,(v,y,w')) & \longmapsto & (w,(v,y,w'))
\end{array} \right)$$
\noindent
where $h(w) = s(m), v
\in \psi_m, y \in \nu_{M}, w' \in -\phi_m, w \in N \subset \phi$ and $m
\in M$.  
\end{Lemma}

{\bf Proof of Lemma \ref{obs:fhetod1nm}}
Let $s_t$ be the path of sections $s_t(m) = t \, s(m)$ from $s$ to the
zero
section and let $s_{[0,1]}$ be the corresponding section of $\pi_{\psi}
\times
id_{[0, 1]}$.  Deform 
$h \times id: \phi \times [0, 1] \rightarrow \psi \times [0, 1]$ to 
$H$ so that $H$ is 
transverse to $s_{[0,1]}(M \times [0,1])$ and so that $H|_{\phi \times
(1-\epsilon,1]} =  (h \times id)|_{\phi \times (1-\epsilon,1]}$. 
Then $H^{-1}(s_{[0,1]}(M \times [0, 1]))$ 
is a normal cobordism between $H^{-1}(s_1(B))$ and
$H^{-1}(s_0(B))$.  The latter is the usual definition of the
normal invariant associated to $h$ (see \cite[Theorem 2.23]{MM}) while the
former is the
one we give above.  To check that the normal data are correct one
observes that $\nu_N = \nu_{N \hookrightarrow \phi} \oplus
\pi_{\phi}|_N^*(\nu_{\phi})$.  The first of these bundles is
$(\pi_{\phi}|_N)^*\psi$ by transversality and the fact that $\psi$ is the
normal
bundle to any section of $\psi$.  The second bundle is computed by
\ref{fact}. All bundle maps are the canonical maps associated to 
pullbacks.  $\qed$

\vspace{0.3cm}

We now complete the proof of Lemma \ref{normalinvt}.
Consider the fiber homotopy equivalence $\pi_{\xi}^*f':
\pi_{\xi}^*(\zeta)
\rightarrow  \pi_{\xi}^*(\xi)$ and the diagonal section
$s:S(\xi) \hookrightarrow \pi_{\xi}^*(\xi), s(x)=(x,x)$.  
In general, given a map of manifolds $F:N \rightarrow M$, $F\times id_M$
is
transverse to the diagonal in $M \times M$ with inverse image the graph of
$F$.  
We have a fiberwise version of this situation so that $\pi_{\xi}^*f'$ 
is transverse to $s(S(\xi))$ with inverse image
$N:=\pi_{\xi}^*f'^{-1}(s(S(\xi))) = 
\{(x,z) | x =f(z) \} \subset \pi_{\xi}^*\zeta$. 
As $f'$ is length preserving $|x|=|z|=1$ for all $(x,z) \in N$.  Hence we
identify $S(\zeta)$ and $N$ via the embedding $i:S(\zeta) \hookrightarrow
\pi_{\xi}^*\zeta,i(z)=(f(z),z).$  The bundle projection
$\pi_{\xi}^*(\zeta) \rightarrow S(\xi)$ restricted to $N$ takes pairs
$(f(z),z)$ to $f(z)$ and we therefore identify it with $f:S(\zeta)
\rightarrow
S(\xi)$.  Applying  observation \ref{obs:fhetod1nm} and the fact 
that $g^{*} \, \pi_{\zeta}^{*} = \pi_{\xi}^{*}$ we obtain a degree 
one normal map with domain 
$$\begin{array}{lcl}
\nu_{S(\zeta)} & = & f^*(\pi_{\xi}^*\xi \oplus \nu_{S(\xi)} \oplus 
- \pi_{\xi}^*\zeta) \\ 
& = & f^*\pi_{\xi}^*(\xi \oplus \nu_{S(\xi)} \oplus -\xi \oplus - \zeta)
\\
& = & \pi_{\zeta}^*(\nu_{S(\xi)} \oplus -\zeta) \oplus \pi_{\zeta}^*(\xi
\oplus
-\xi). \\
\end{array}$$
The degree one normal map itself is given as follows.
\begin{equation}
\label{eq:nof'}
\eta^O(\{\pi_{\xi}^*f'\}) = \left( \begin{array}{ccc}
\pi_{\zeta}^*(\nu_{S(\xi)} \oplus - \zeta \oplus  \xi \oplus - \xi ) &
\longrightarrow 
& \pi_{\xi}^*(\nu_{S(\xi)} \oplus - \zeta  \oplus \xi \oplus - \xi ) \\
\downarrow & & \downarrow  \\
S(\zeta) & \longrightarrow & S(\xi) \\
(z,(y,z',x,x')) & \longmapsto & (f(z),(y,z',x,x'))
\end{array} \right)
\end{equation}
\noindent
Since the factor $\xi \oplus - \xi$ is a trivial summand onto which we map
with the identity, we can ignore it.  Comparing equalities 
(\ref{eq:nof}) and
(\ref{eq:nof'}), we
see  that $\eta^O(\{\pi_{\xi}^*f'\}) = \eta^O(\{ f\})$.  Finally note 
that the definition of $f'$ implies that $\eta^O(\{\pi_{\xi}^*f'\}) = 
\pi_{\xi}^*([f])\,.$
$\qed$

\begin{Remark}
In our case the base space of the sphere bundles is four dimensional.  
Note that we do not use any surgery theory on this four dimensional 
base space.  We only describe the set of normal invariants.
\end{Remark}

\begin{Corollary}
\label{normalinvts}
The fiber homotopy equivalences $f_j:M_{m+12j,n} \rightarrow M_{m,n}$ have
normal invariant $j \in \NO(M_{m,n}) \cong \Z_n.$
\end{Corollary}

{\bf Proof.}  Note that
$\NO(S^4) = \pi_4(G/O)$ and the homotopy exact sequence for the fibration
$G/O ~\hookrightarrow~BO~\rightarrow~BG\,$ at $\,\pi_4\,$ is 
$0 \rightarrow \Z \rightarrow \Z \rightarrow \Z_{24} \rightarrow 0.$
Thus the value of $[f_j]$ in $[S^4,G/O]$ may be computed simply by taking
the difference
of the stable bundles $\xi_{m+12j,n} \oplus \epsilon$ and
$\xi_{m,n} \oplus
\epsilon$ in $\pi_4(BO)$ and then dividing by $24$.  But one computes that
this is
precisely $(n+2(m+12j) - (n+2m))/24 = j$.  The map $\pi_{\xi}^*: \Z \cong
[S^4,G/O]
\rightarrow [M_{m,n},G/O] \cong  \Z_n$ is identical with the map on
cohomology
$\pi_{\xi}^* : H^4(S^4) \rightarrow H^4(M_{m,n})$ and thus
$\pi_{\xi}^*[f_j] = j
\in \Z_n = [M_{m,n},G/O] = \NO(M_{m,n})$.
$\qed$

\medskip

{\bf Proof of Theorem \ref{hom}}
Let $\Theta^7$ denote the group of exotic seven
spheres.  Since $\mu$ is additive with respect to connected sums and 
distinguishes all homotopy
$7$-spheres, the action of $\Theta^7$ on ${\mathcal N}^O(M_{m,n})$, 
$\,\Sigma \cdot [f:M
\rightarrow M_{m,n}] = [f:M \,\sharp \,\Sigma \rightarrow
M_{m,n}]$,  is free.  The quotient under this action defines the map
$\eta^O$
to the smooth normal invariants.
$$
0 \longrightarrow \Theta^7 \longrightarrow {\mathcal S}^{O}(M_{m,n})
\stackrel{\eta^O}{\longrightarrow} {\mathcal N}^{O}(M_{m,n})
\longrightarrow 0.
$$
We see that all normal
invariants of $M_{m,n}$ are normally bordant to homotopy equivalences 
and distinct elements of the
structure set with the same normal invariants differ only by a 
connected sum with an exotic $7$-sphere in the domain.  But we have 
constructed the normal invariant set in Lemma \ref{normalinvt}.  It
follows 
that a given seven dimensional manifold $M$ is homotopy 
equivalent to $M_{m,n}$ if and only if
$M$ is homeomorphic to $M_{m+12j,n}$ for some $j$.  Moreover, in Theorem
\ref{topcor} we  gave conditions on $m$ and $m'$ for a classification up
to almost diffeomorphism.  It remains to check that applying Theorem
\ref{top}
yields the relation given in Theorem \ref{hom} and we do this now.

\vspace{0.3cm}

\noindent \emph{Case $n =2^0\,q\,,$ $q$ odd and $\halfp$ classifies.}
$$\begin{array}{lcl}
M_{m',n} \simeq M_{m,n} & \Leftrightarrow & M_{m',n} \cong
M_{m+12\,j,n}~\mbox{for some} ~j. \\
& \Leftrightarrow & \alpha ~ m' \equiv m
+12\,j~\mod n~\mbox{for some}~j~\mbox{and some}~\alpha \in
A^+(n). \\
& \Leftrightarrow & \alpha ~ m' \equiv  m~\mod (n,12)~
\mbox{for some}~\alpha \in A^+(n). \\
& \Leftrightarrow & \alpha ~ m' \equiv  m~\mod (n,12)~\mbox{for
some}~\alpha \in
A^+((n,12)).
\end{array}$$

\noindent
\emph{Case $n = 2^aq$ where $a=2$ and $m$ odd or $a>2$ and $m$ even, $\halfp$
classifies and
$M_{m+\epsilon \,\frac{n}{2},n} \cong M_{m,n}$ for $\epsilon =0,1$}.
$$\begin{array}{lcl}
M_{m',n} \simeq M_{m,n} & \Leftrightarrow & M_{m',n}
\cong M_{m+12\,j+\epsilon\,\frac{n}{2},n}~\mbox{for 
some}~j ~\mbox{and for}~\epsilon=0,1. \\ 
&\Leftrightarrow & \alpha ~ m' \equiv
m+12\,j+\epsilon \frac{n}{2}~\mod \frac{n}{2} ~\mbox{for
some}~j,\epsilon=0, 1~
\mbox{and}~\alpha \in A^+(n). \\
&\Leftrightarrow & \alpha ~ m' \equiv m +12\,j~\mod n~\mbox{for
some}~\alpha \in
A^+(n).\\
&\Leftrightarrow & \alpha ~ m' \equiv m~\mod (n,12)~\mbox{for some}~\alpha
\in
A^+((n,12)).
\end{array}$$

\noindent \emph{Case $n =2^a \, q, a=1$ or $a=2$ and $m$ even or $a>2$ 
and $m$ odd.}
$$\begin{array}{lcl}
M_{m',n} \simeq M_{m,n} & \Leftrightarrow & M_{m',n} \cong
M_{m+12\,j,n}~\mbox{for some} ~j. \\
& \Leftrightarrow & \alpha \, m' \equiv m +12\,j \,\mod n ~\mbox{for 
some}\,j\,\mbox{and} \,\alpha \in A^+(n),
\alpha \equiv \pm 1 \,\mod 2^a. \\
& \Rightarrow & \alpha ~ m' \equiv m~\mod (n,12)~\mbox{for some}~\alpha~\in
A^+(n).\\
& \Leftrightarrow & \alpha ~ m' \equiv m~\mod (n,12)~\mbox{for
some}~\alpha \in
A^+((n,12)).
\end{array}$$
\noindent
To reverse the third implication we must show that if $\alpha ~ m' \equiv
m +12\,j~\mod n$ for some
$\alpha \in A^+(n)$ then $\bar{\alpha}$ and $\bar{j}$ can be chosen 
so that $\bar{\alpha} ~
m' \equiv m + 12\,\bar{j}~\mod n$ and $\bar{\alpha} \equiv \pm
1~\mod 2^a$.  This is
evident for $a=1,2$ so we assume that $a>2$ and $\alpha \not\equiv \pm
1~\mod n$.  But
now $\frac{n}{2} \equiv 12\,k~\mod n$ for some $k$ and, since $m$ is odd,
$(\alpha+\frac{n}{2}) m'
\equiv m + 12\,(j+k)~\mod n$.  Moreover $(\alpha +\frac{n}{2})^2 \equiv
1~\mod n$, hence
$\bar{\alpha} = \alpha + \frac{n}{2}$ and $\bar{j}  = j + k$ suffice.

\vspace{0.3cm}
\noindent 
In the orientation reversing case, $M_{m',n} \simeq -M_{m,n}$ if and only
if 
$M_{m',n} \cong -M_{m+12\,j,n}$ for some $j$.  But we know from Theorem
\ref{top} that this
only happens when $n=2^{\epsilon}p_1^{i_1} \dots p_k^{i_k}$,
\noindent
$~p_i~\equiv~1~\mod 4$ and $\epsilon = 0,1$.  When
$\epsilon = 0$ there is only one oriented homotopy type and it admits an
orientation reversing homotopy equivalence by Theorem \ref{top}.  When
$\epsilon=1$ there are two oriented homotopy types and these are
orientation
reversing homotopy  equivalent again by Theorem \ref{top}.  This concludes
the
proof of Theorem \ref{hom}.

\section{Remarks on the classification}
\subsection{ The integer 4 is not favorable}
The significance of orientation in the classification is strikingly
illustrated by the pair of manifolds $M_{-1,2}$ and
$M_{0,2}$.  The first is the total space of the unit tangent bundle of 
$S^4$ and the second is the total space of 
``twice'' the Hopf fibration.  It follows from Theorem \ref{topcor} that
these
bundles are orientation reversing diffeomorphic, a fact well known to
geometers, see for example \cite[p.192, example g]{Ri}, but by Theorem
\ref{hom} they are not 
orientation preserving
homotopy equivalent.  Notice this situation implies
that the oriented homotopy types represented by $M_{-1,2}$ and $M_{0,2}$
do
not admit orientation reversing homotopy equivalences. If they did then
they
would be orientation preserving homotopy equivalent.  James and Whitehead
(\cite[p.151]{JW2}) called an integer $k$ favorable when the existence of
an
orientation reversing homotopy equivalence between $(k-1)$-sphere bundles
over $S^k$ with Euler number $2$ implies the existence of an orientation
preserving homotopy equivalence between the bundles.  We have thus shown 
that \emph {$4$ is not favorable.}

\subsection{The total spaces $M_{m,n}$ versus the pairs $(M_{m,n},S^3)$}
 We noted above that the homotopy classification
of the total spaces and the homotopy classification of the pairs differ if
and only
if
$(n,12)=12$.   For example, $M_{1,12}$ and $M_{5,12}$ are homeomorphic
but
the pairs $(M_{1,12},S^3)$ and $(M_{5,12},S^3)$ are not homotopy
equivalent. 
Algebraically, the homeomorphism takes a fiber to five times a fiber and
thus
pays no respect to the bundle structures of domain and range.  In
\cite[Theorem 3.7]{JW0} James and Whitehead give connectivity conditions
under
which maps of bundle-like spaces are homotopic to maps of bundle fiber
pairs. 
This example shows that their results are sharp.


\begin{thebibliography}{9999999}
	
\bibitem[D] {D} M. Davis, {\em Some group actions on homotopy spheres
	of dimension seven and fifteen}, Amer. J. of Math. {\bf 104}
	(1982) 59-90.
\bibitem[EK] {EK} J. Eells and N. Kuiper, {\em An invariant for
	certain smooth manifolds}, Annali di Math. {\bf 60} (1962)
	93-110.
\bibitem[GM] {GM} D. Gromoll and W. Meyer, {\em An exotic sphere with
	nonnegative sectional curvature}, Ann. of Math. {\bf 100}
	(1974) 401-406.
\bibitem[GZ] {GZ} K. Grove and W. Ziller, {\em Curvature and symmetry
	of Milnor spheres}, to appear.
\bibitem[H]{H} A.~Hatcher, {\em A proof of the Smale Conjecture,
$Diff(S^3) \simeq O(4)$}, Ann. of Math. {\bf 117} (1983) 553-607.
\bibitem[JW0] {JW0} I. M. James and J. H. C. Whitehead, {\em Note on Fiber
Spaces}, Proc. London Math. Soc. {\bf 3} (1954) 129-137.
\bibitem[JW1] {JW1} I. M. James and J. H. C. Whitehead, {\em The
homotopy theory of spheres bundles over spheres (I)},
Proc. London Math. Soc. {\bf 3} (1954) 196-218.
\bibitem[JW2] {JW2} I. M. James and J. H. C. Whitehead, {\em The
	homotopy theory of spheres bundles over spheres (II)},
	Proc. London Math. Soc. {\bf 5} (1955) 148-166.
\bibitem[K]{K} M.~Kreck, {\em Surgery and Duality}, Ann. of Math. {\bf
149}
no.3 (1999) 707-754.
\bibitem[KiSh]{KiSh} 
N.~Kitchloo and K.~Shankar, {\em On complexes equivalent to
$S^3$-bundles over $S^4$}, preprint.
\bibitem[KS] {KS} M. Kreck and S. Stolz, {\em Some nondiffeomorphic 
homeomorphic homogeneous $7$-manifolds with positive sectional 
curvature}, J. Diff. Geom. {\bf 33} (1991) 465-486.
\bibitem[L] {L} W. LeVeque, {\em Fundamentals of number theory}, 
Addison-Wesley Publishing Comp. (1977).
\bibitem[M] {M} J. Milnor, {\em On manifolds homeomorphic to the
	$7$-sphere}, Annals of Math. {\bf 64} (1956) 399-405.
\bibitem[Ma] {Ma} M.Masuda, {\em Smooth group actions on sphere bundles 
over spheres and on Brieskorn manifolds}, Proc. of Amer. Math. Soc.
{\bf 92} (1984) 119-124.
\bibitem[MM] {MM} I. Madsen and J. Milgram, {\em The classifying
	spaces for surgery and cobordism of manifolds}, Annals of
	Mathematics Studies, Princeton University Press (1979).
\bibitem[S] {S} S. Sasao, {\em On homotopy type of certain complexes}, 
Topology {\bf 3} (1965) 97-102. 
\bibitem[Ri] {Ri} A. Rigas, {\em Some bundles of non-negative 
curvature}, Math. Ann. {\bf 232} (1978) 187-193.  
\bibitem[T] {T} I. Tamura, {\em Homeomorphism classification of total
spaces of sphere bundles over spheres}, J. of Math. Soc. Japan {\bf 10} 
(1958) 29-43. 
\bibitem[Wa1]{Wa1} C. T. C. Wall, {\em Classification problems in
differential topology - IV}, Topology, {\bf 6} (1967) 273-296.
\bibitem[Wa2]{Wa2} C. T. C. Wall {\em Quadratic forms on finite groups and
related topics}, Topology, {\bf 2} 91964) 281-298.
\bibitem[Wh]{Wh} G.~W.~Whitehead, {\em Elements of homotopy theory}, 
Graduate Texts in Mathematics, Springer-Verlag (1978).
\bibitem[Wi] {Wi} D. Wilkens,  {\em Closed $(s-1)$-connected 
 $(2s+1)$-manifolds, $s = 3,7$},
 Bull. London Math. Soc. {\bf 4} (1972)  27-31.  
\end{thebibliography}
\end{document}